\documentclass[12pt]{article}
\usepackage{graphicx}
\usepackage{amsmath}
\usepackage{amsfonts}
\usepackage{amssymb}
\newtheorem{theorem}{Theorem}

\newtheorem{corollary}[theorem]{Corollary}

\newenvironment{proof}[1][Proof]{\textbf{#1.} }{\ \rule{0.5em}{0.5em}}

\def\text{\hbox} 

\def\a{\alpha}

\def\f{\phi}

\def\o{\omega}

\def\F{\Phi}

\def\S{{\bf S}}

\def\A{{\cal A}}
\def\B{{\cal B}}
\def\A'{{\cal A}'}
\def\B'{{\cal B}'}

\def\Z{{\mathbf Z}}

\begin{document}

\title{Cyclic presentations of groups and cyclic branched coverings of
$(1,1)$-knots \footnote{Work performed under the auspices of the
G.N.S.A.G.A. of I.N.d.A.M. (Italy) and the University of Bologna,
funds for selected research topics.}}
\author{Michele Mulazzani}

\maketitle
\begin{abstract}
{In this paper we study the connections between cyclic
presentations of groups and cyclic branched coverings of
$(1,1)$-knots. In particular, we prove that every $n$-fold
strongly-cyclic branched covering of a $(1,1)$-knot admits a
cyclic presentation for the fundamental group encoded by a
Heegaard
diagram of genus $n$. }

\end{abstract}


\section{Introduction and preliminaries}

The problem of determining whether a balanced presentation of a
group is geometric (i.e. induced by a Heegaard diagram of a closed
orientable 3-manifold) is of considerable interest in geometric
topology and has already been examined by many authors (see
\cite{Gr}, \cite{Mo}, \cite{Ne}, \cite{OS1}, \cite{OS2},
\cite{OS3}, \cite{St}). Furthermore, the connections between
cyclic coverings of $\S^3$ branched over knots and cyclic
presentations of groups induced by suitable Heegaard diagrams have
recently been discussed in several papers (see \cite{BKM},
\cite{CHK}, \cite{Du}, \cite{HKM1}, \cite{HKM2}, \cite{Ki},
\cite{KKV1}, \cite{KKV2}, \cite{MR}, \cite{VK}).

Note that a finite balanced presentation of a group
$<x_1,\ldots,x_n\mid r_1,\ldots,r_n>$ is said to be a {\it cyclic
presentation\/} if there exists a word $w$ in the free group $F_n$
generated by $x_1,\ldots,x_n$ such that the relators of the
presentation are $r_k=\theta_n^{k-1}(w)$, $k=1,\ldots,n$, where
$\theta_n :F_n\to F_n$ denotes the automorphism defined by
$\theta_n (x_i)=x_{i+1}$ (mod $n$), $i=1,\ldots,n$. This cyclic
presentation (and the related group) will be denoted by $G_n(w)$.
For further details see \cite{Jo}.

\medskip

We list the most interesting examples:
\begin{itemize}
\item the Fibonacci group
$F(2n)=G_{2n}(x_1x_2x_3^{-1})=G_n(x_1^{-1}x_2^2x_3^{-1}x_2)$ is
the fundamental group of the $n$-fold cyclic covering of $\S^3$
branched over the figure-eight knot, for all $n>1$ (see
\cite{HKM2});
\item the Sieradsky group
$S(n)=G_n(x_1x_3x_2^{-1})$ is the fundamental group of the
$n$-fold cyclic covering of $\S^3$ branched over the trefoil knot,
for all $n>1$ (see \cite{CHK});
\item the fractional Fibonacci group
$\widetilde F_{l,k}(n)=G_n((x_1^{-l}x_2^l)^kx_2(x_3^{-l}x_2^l)^k)$
is the fundamental group of the $n$-fold cyclic covering of $\S^3$
branched over the genus one two-bridge knot with Conway
coefficients $[2\,l,-2k]$, for all $n>1$ and $l,k>0$ (see
\cite{VK}).
\end{itemize}
Moreover, all the above cyclic presentations are geometric (i.e.,
they arise from suitable Heegaard diagrams).

In order to investigate these relations, Dunwoody introduced in
\cite{Du} a class of Heegaard diagrams depending on six integers,
having cyclic symmetry and defining cyclic presentations for the
fundamental group of the represented manifolds. In \cite{GM} it
has been shown that the 3-manifolds represented by these diagrams
are cyclic coverings of lens spaces, branched over $(1,1)$-knots
(also called genus one 1-bridge knots). As a corollary, it has
been proved that for some determined cases the manifolds turn out
to be cyclic coverings of $\S^3$, branched over suitable knots.
This gives a positive answer to a conjecture formulated by
Dunwoody in \cite{Du}. Section 2 resumes the main statements of
\cite{GM} concerning this topic.

The above results suggest that cyclic presentations of groups are
actually related to cyclic branched coverings of $(1,1)$-knots. As
a basic result in this direction, we prove in Section 3 that every
$n$-fold strongly-cyclic branched covering of a $(1,1)$-knot
admits a Heegaard diagram of genus $\,n$ which encodes a cyclic
presentation for the fundamental group. The definition of
strongly-cyclic branched coverings of $(1,1)$-knots will be
introduced in Section 3. It is interesting to note that the
construction used to prove our main result is strictly related to
the concept of $p$-symmetric Heegaard splittings introduced by
Birman and Hilden in \cite{BH}.

In what follows, we shall deal with $(1,1)$-knots, i.e. knots in
lens spaces (possibly in $\S^3$), which admit a certain
decomposition. A knot $K$ in a lens space $L(p,q)$ is called a
$(1,1)$-{\it knot\/} (or also a {\it genus one 1-bridge knot\/}) if
there exists a Heegaard splitting of genus one
$(L(p,q),K)=(T,A)\cup_{\f}(T',A')$, where $T$ and $T'$ are solid
tori, $A\subset T$ and $A'\subset T'$ are properly embedded
trivial arcs, and $\f:(\partial T,\partial A)\to(\partial
T',\partial A')$ is the attaching homeomorphism. This means that there
exists a disk $D\subset T$ (resp. $D'\subset T'$) with $A\cap
D=A\cap\partial D=A$ and $\partial D-A\subset\partial T$ (resp.
$A'\cap D'=A'\cap\partial D'=A'$ and $\partial
D'-A'\subset\partial T'$).

\bigskip

\begin{figure}[ht]
 \begin{center}
 \includegraphics*[totalheight=3cm]{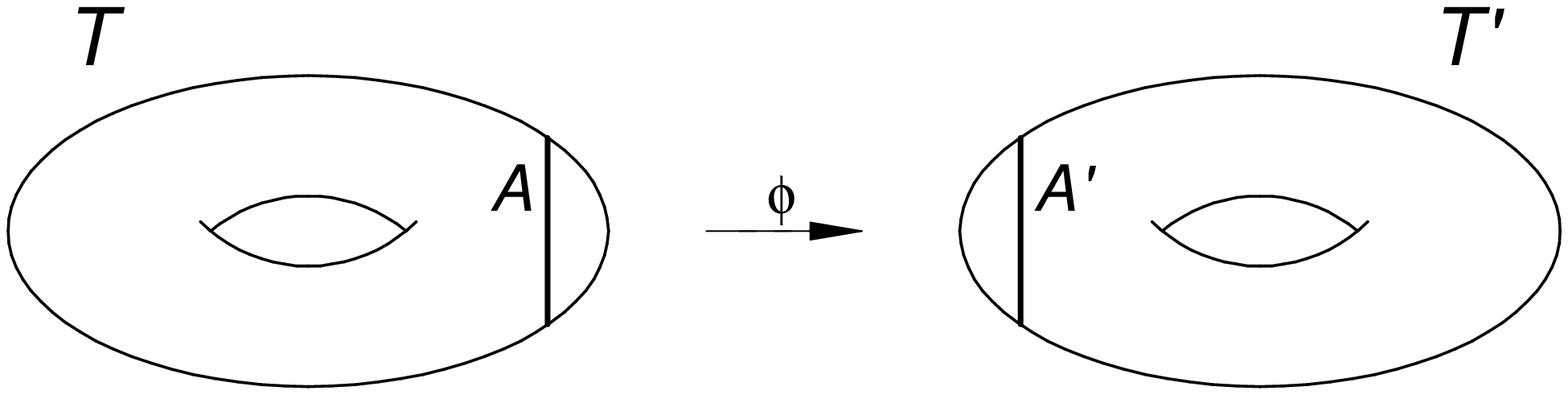}
 \end{center}
 \caption{A $(1,1)$-knot decomposition.}

 \label{Fig. 3}

\end{figure}

Notice that $(1,1)$-knots are only a particular case of the notion
of $(g,b)$-links in closed orientable 3-manifolds (see \cite{Do}
and \cite{Ha}), which generalize the classical concept of bridge
decomposition of links in $\S^3$.

The class of $(1,1)$-knots has recently been investigated in
several papers (see \cite{Be, Fu, Ga, Ga2, GM, Ha, Ha2, MS, MSY,
MSY2, W1, Wu, W2}) and appears very important in the light of some
results and conjectures involving Dehn surgery on knots. It is
well known that the subclass of $(1,1)$-knots in $\S^3$ contains
all torus knots (trivially) and all 2-bridge knots (i.e.,
$(0,2)$-knots) \cite{MS}.

\section{Dunwoody manifolds}

M.J. Dunwoody introduced in \cite{Du} a class of Heegaard diagrams
having a cyclic symmetry, depending on six integers $a,b,c,n,r,s$
such that $n>0$, $a,b,c\ge 0$ and $a+b+c>0$ (see Figure 2). This
construction gives rise to a wide class of closed orientable
3-manifolds $D(a,b,c,n,r,s)$, called {\it Dunwoody manifolds\/},
admitting geometric cyclic presentations for their fundamental
groups.

The diagram is an Heegaard diagram of genus $n$. It contains $n$
upper cycles $C'_1,\ldots,C'_n$ and $n$ lower cycles
$C''_1,\ldots,C''_n$, each having $d=2a+b+c$ vertices. For each
$i=1,\ldots,n$, the cycle $C'_i$ (resp. $C''_i$) is connected to
the cycle $C'_{i+1}$ (resp. $C''_{i+1}$) by $a$ parallel arcs, to
the cycle $C''_{i}$ by $c$ parallel arcs and to the cycle
$C''_{i+1}$ by $b$ parallel arcs. The cycle $C'_i$ is glued to the
cycle $C''_{i-s}$ (mod $n$) so that equally labelled vertices are
identified together (the labelling of the cycles is pointed out in
Figure 3).

It is evident that the diagram (as well as the identification
rule) is invariant with respect to a cyclic action of order $n$.

\begin{figure}[ht]
 \begin{center}
 \includegraphics*[totalheight=7cm]{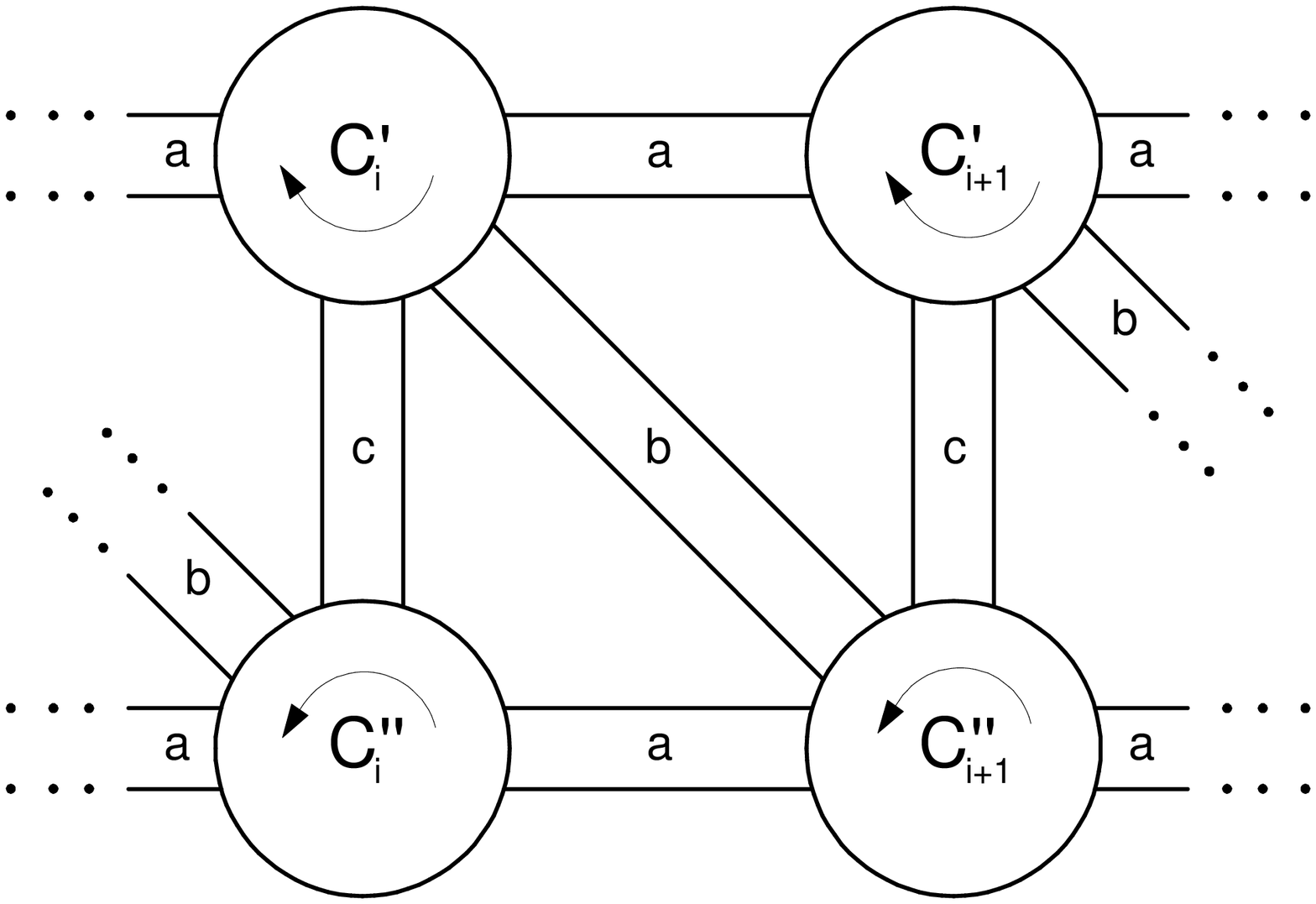}
 \end{center}
 \caption{Heegaard diagram of Dunwoody type.}

 \label{Fig. 1}

\end{figure}


\begin{figure}[ht]
 \begin{center}
 \includegraphics*[totalheight=8cm]{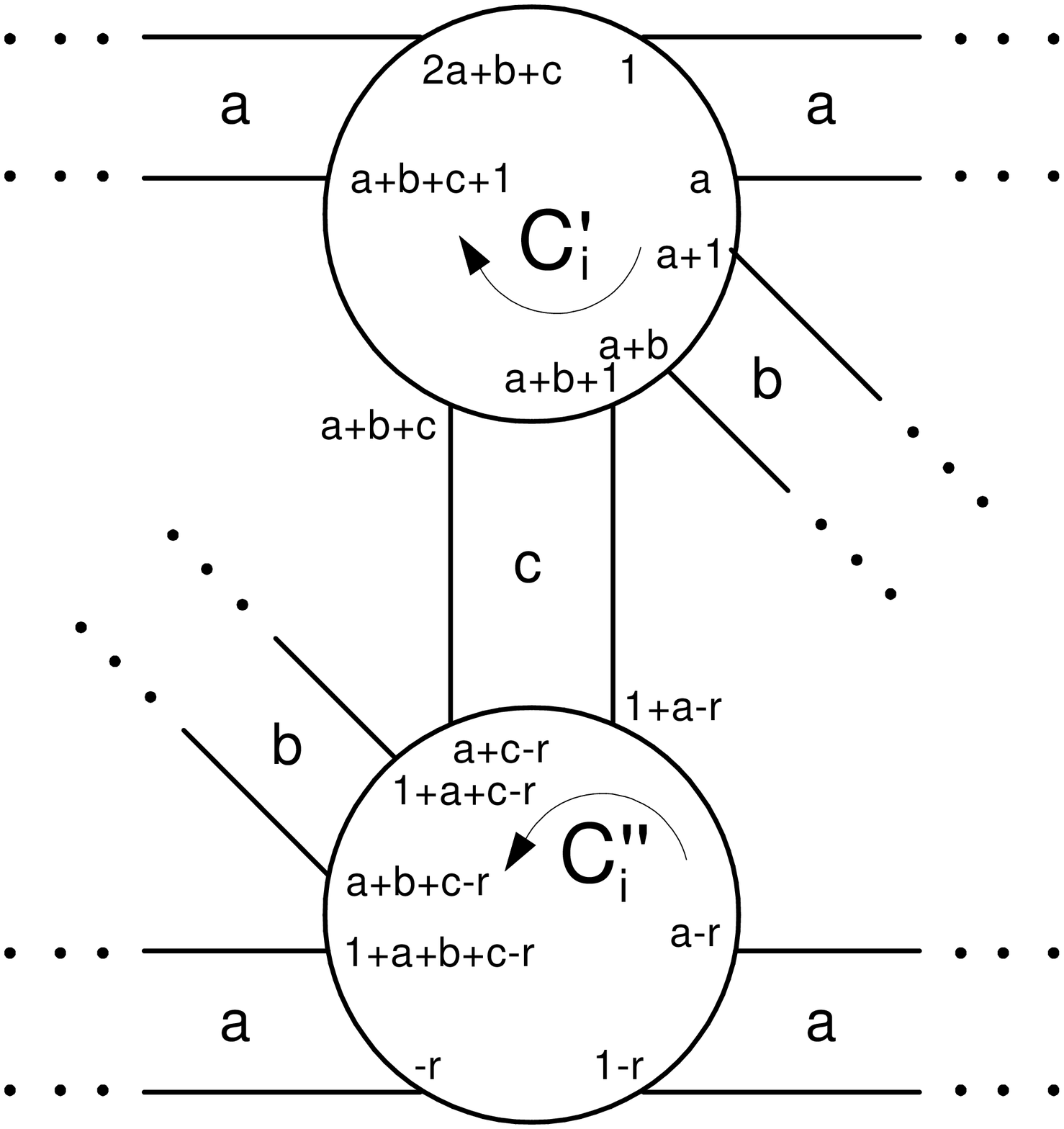}
 \end{center}
 \caption{}

 \label{Fig. 2}

\end{figure}

In \cite{GM} it has been shown that
each Dunwoody manifold is a cyclic covering of a lens space
(possibly $\S^3$), branched over a $(1,1)$-knot.

\begin{theorem} \label{Theorem 3} \cite{GM}
The Dunwoody manifold $D(a,b,c,n,r,s)$ is the $n$-fold cyclic
covering of a lens space $D'$ (possibly $\S^3$) branched
over a $(1,1)$-knot $K\subset D'$, both only depending on the
integers $a,b,c,r$.
\end{theorem}

As a consequence, for certain particular values of the parameters
a Dunwoody manifold turn out to be a cyclic covering of $\S^3$
branched over a knot. This gives a positive answer to a conjecture
formulated by Dunwoody in \cite{Du}, which has also been
independently proved in \cite{SK}.

\begin{corollary} \label{Corollary 4} \cite{GM} Each
Dunwoody manifold $D(a,b,c,n,r,s)$ of Section 3 of \cite{Du} is an
$n$-fold cyclic covering of $\bf S^3$, branched over a
$(1,1)$-knot $K\subset \S^3$, which only depends on $a,b,c,r$.
\end{corollary}

\section{Strongly-cyclic branched coverings of $(1,1)$-knots}

As well known, an $n$-fold cyclic branched covering between two
orientable closed manifolds $f:M\to N$, with branching set $L$, is
completely defined by an epimorphism $\o_{f}:H_1(N-L)\to\Z_n$,
where $\Z_n$ is the cyclic group of order $n$. If $N=\S^3$ and the
branching set is a knot $K$, the covering is uniquely determined,
up to covering equivalence, since $H_1(\S^3-K)\cong\Z$ and the
homology class $[m]$ of a meridian loop around the knot have to be
mapped by $\o_{f}$ in a generator of $\Z_n$. Therefore the index
of the branching set is exactly $n$.

Obviously, this property does not hold for a $(1,1)$-knot in a
lens space. Moreover, we would like to obtain cyclic branched
coverings producing a cyclic presentation for the fundamental
group of the manifold. In order to achieve this, we will select
cyclic branched coverings of ``special type'', and this will be a
very natural generalization of the case of knots in $\S^3$.

An $n$-fold cyclic covering of $L(p,q)$ branched over the
$(1,1)$-knot $K$ will be called {\it strongly-cyclic\/} if the
branching index of $K$ is $n$. This means that the homology class
of a meridian loop $m$ around $K$ is mapped by $\o_{f}$ in a
generator of $\Z_n$ (up to covering equivalence we can always
suppose $\o_{f}[m]=1$).

Strongly-cyclic branched coverings of $(1,1)$-knots appear to be a
suitable tool for producing 3-manifolds with fundamental group
admitting cyclic presentation. For example, it is easy to see that
all Dunwoody manifolds are coverings of this type.

\begin{theorem} \label{Main Theorem}
Every $n$-fold strongly-cyclic branched covering of a $(1,1)$-knot
admits a Heegaard diagram of genus $n$, which induces a cyclic
presentation of the fundamental group of the manifold.
\end{theorem}
\begin{proof}
Let $f:(M,f^{-1}(K))\to (L(p,q),K)=(T,A)\cup_{\f}(T',A')$ be an
$n$-fold strongly-cyclic branched covering of the $(1,1)$-knot
$K$. Then $Y_n=f^{-1}(T)$ and $Y_n'=f^{-1}(T')$ are both
handlebodies of genus $n$. Moreover, $f^{-1}(A)$ and $f^{-1}(A')$
are both properly embedded arcs in $Y_n$ and $Y_n'$ respectively.
We get a genus $n$ Heegaard splitting
$(M,f^{-1}(K))=(Y_n,f^{-1}(A))\cup_{\F}(Y'_n,f^{-1}(A'))$, where
$\F:\partial Y_n\to\partial Y_n'$ is the lifting of $\f$ with
respect to $f$. Let $m\subset T-A$ be a meridian loop around $A$
and $\a\subset T-A$ be a generator of $\pi_1(T,P)$ such that
$\o_{f}[\a]=0$, where the base point $P$ is any point of $m$. The
loop $\a$ exists: take a generator $\a'\subset T-A$ of
$\pi_1(T,P)$; if $\o_{f}[\a']=k$ then choose any $\a$ homotopic to
$\a'm^{-k}$. Moreover, take a point $Q\in A$ and let $\gamma$ be
an arc from $P$ to $Q$ such that $\gamma\cap A=Q$. Then
$f^{-1}(Q)$ is a single point $*\in f^{-1}(A)$ and $f^{-1}(P)$
consists of $n$ points $P_1,\ldots,P_n$. For $i=1,\ldots,n$, let
$\tilde\a_i$ and $\tilde\gamma_i$ be the lifting (with respect to
$f$) of $\a$ and $\gamma$ respectively, both containing $P_i$.
Then the $n$ loops
$\bar\a_1=\tilde\gamma_1^{-1}\cdot\tilde\a_1\cdot\tilde\gamma_1,
\ldots,\bar\a_n=\tilde\gamma_n^{-1}\cdot\tilde\a_n\cdot\tilde\gamma_n$
generate $\pi_1(Y_n,*)$ and they are cyclically permutated by a
generator $\Psi$ of the group of covering transformations. Let
$E'$ be a meridian disk for the torus $T'$ such that $E'\cap
A'=\emptyset$, then $f^{-1}(E')$ is a system of meridian disks
$\{\tilde E_1',\ldots,\tilde E_n'\}$ for the handlebody $Y_n'$,
and they are cyclically permutated by $\Psi$. The curves
$\F^{-1}(\partial\tilde E_1'),\ldots,\F^{-1}(\partial\tilde E_n')$
give the relators for the presentation of $\pi_1(M,*)$ induced by
the Heegaard splitting. Since both generator and relator curves
are cyclically permuted by $\Psi$, we get the statement.
\end{proof}

\medskip

Obviously several problems arise:
\begin{itemize}
\item study the relations between the attaching homeomorphism $\,\f$
producing $K$ and the monodromy $\o_{f}$ of the strongly-cyclic
branched covering;
\item find the word $w$ associated to the cyclic presentation,
starting from the $(1,1)$-knot description;
\item find some strongly-cyclic
branched covering of a $(1,1)$-knot (possibly in $\S^3$) which is
not a Dunwoody manifold.
\end{itemize}

A discussion of the first two problems can be found in \cite{CM}.

\vspace{15 pt} {MICHELE MULAZZANI, Department of Mathematics,
University of Bologna, I-40127 Bologna, ITALY, and C.I.R.A.M.,
Bologna, ITALY. E-mail: mulazza@dm.unibo.it}

\end{document}